\documentclass[12pt]{amsart}
\usepackage{amsmath,amscd,amssymb,amsfonts}
\newcommand{\h}{\hbox}
\newcommand{\q}{\quad}
\newcommand{\nin}{\noindent}
\newcommand{\bs}{\par\bigskip}
\newcommand{\ms}{\par\medskip}
\newcommand{\sk}{\par\smallskip}

\newcommand{\msum}{\h{$\sum$}}
\newcommand{\mopl}{\h{$\bigoplus$}}
\newcommand{\mcap}{\h{$\bigcap$}}
\newcommand{\mcup}{\h{$\bigcup$}}
\newcommand{\mprod}{\h{$\prod$}}
\newcommand{\mcoprod}{\h{$\coprod$}}
\newcommand{\ssb}{\raise.15ex\h{${\scriptscriptstyle\bullet}$}}
\newcommand{\ssc}{\,\raise.15ex\h{${\scriptstyle\circ}$}\,}
\newcommand{\C}{{\mathbf C}}
\newcommand{\D}{{\mathbf D}}
\newcommand{\N}{{\mathbf N}}
\newcommand{\PP}{{\mathbf P}}
\newcommand{\Q}{{\mathbf Q}}
\newcommand{\R}{{\mathbf R}}
\newcommand{\Z}{{\mathbf Z}}

\newcommand{\CC}{{\mathcal C}}
\newcommand{\DD}{{\mathcal D}}
\newcommand{\HH}{{\mathcal H}}
\newcommand{\K}{{\mathcal K}}
\newcommand{\LL}{{\mathcal L}}
\newcommand{\M}{{\mathcal M}}
\newcommand{\OO}{{\mathcal O}}
\newcommand{\SSS}{{}\,\overline{\!S}{}}

\newcommand{\uu}{\widetilde{u}}
\newcommand{\lc}{{\it loc.~cit.}}
\newcommand{\Gr}{{\rm Gr}}
\newcommand{\DR}{{\rm DR}}
\newcommand{\Ext}{{\rm Ext}}
\newcommand{\bl}{\bigl}
\newcommand{\br}{\bigr}
\newcommand{\simto}{\buildrel\sim\over\longrightarrow}
\newcommand{\into}{\hookrightarrow}
\begin{document}
\title{Some remarks on the semi-positivity theorems}
\author[O. Fujino]{Osamu Fujino}
\address{Department of Mathematics, Faculty of Science, Kyoto University, Kyoto 606-8502, Japan}
\email{fujino@math.kyoto-u.ac.jp}
\author[T. Fujisawa]{Taro Fujisawa}
\address{Tokyo Denki University, School of Engineering, Department of Mathematics, Tokyo, Japan}
\email{fujisawa@mail.dendai.ac.jp}
\author[M. Saito]{Morihiko Saito}
\address{RIMS Kyoto University, Kyoto 606-8502 Japan}
\email{msaito@kurims.kyoto-u.ac.jp}
\begin{abstract}
We show that the dualizing sheaves of reduced simple normal crossings pairs have a canonical weight filtration in a compatible way with the one on the corresponding mixed Hodge modules by calculating the extension classes between the dualizing sheaves of smooth varieties. Using the weight spectral sequence of mixed Hodge modules, we then reduce the semi-positivity theorem for the higher direct images of dualizing sheaves to the smooth case where the assertion is well known. This may be used to simplify some constructions in a recent paper of Y. Kawamata. We also give a simple proof of the semi-positivity theorem for admissible variations of mixed Hodge structure in [FF] by using the theories of Cattani, Kaplan, Schmid, Steenbrink, and Zucker. This generalizes Kawamata's classical result in the pure case.

\end{abstract}
\maketitle
\centerline{\bf Introduction}
\bs\nin
Let $(X,D)$ be a reduced simple normal crossing pair (i.e., $(X,D)$ is Zariski-locally isomorphic to $(X',X'\cap D')$ with $X',D'$ reduced simple normal crossing divisors on a smooth variety $Y$ having no common irreducible components and such that their union is a normal crossing divisor). For $k>0$, set
$$X^{[k]}:=\{x\in X\mid{\rm mult}_xX\ge k\}^{\sim},$$
where $Z^{\sim}$ denotes the normalization of $Z$ for an algebraic variety $Z$. Then $X^{[k]}$ is the disjoint union of the intersections of $k$ irreducible components of $X$, and is smooth. (There is a shift of the index $k$ by one if we compare it with [Kaw3]. There is no shift with the one in [De2], [SZ] although different notation is used there.) We have a reduced simple normal crossing divisor $D^{[k]}\subset X^{[k]}$ defined by the pull-back of $D$ by the natural morphism $X^{[k]}\to X$. For $l\ge 0$, set
$$D^{[k,l]}:=\{x\in X^{[k]}\mid{\rm mult}_xD^{[k]}\ge l\}^{\sim}.$$
Note that $D^{[k,0]}=X^{[k]}$, and our $D^{[k,l]}$ is contained in $X^{[k-1],[l]}$ in [Kaw3] which is much larger in general (and there is a shift of indices as noted above).
\sk
Let $U=X\setminus D$ with $j:U\into X$ the natural inclusion. Set $n=\dim X$. Let $\D\Q_U$ be the dual of the constant sheaf $\Q_U$ in the derived category $D^b_c(X,\Q)$. Since $j$ is an affine immersion, $j_!\Q_U[n]$ and its dual $\R j_*\D\Q_U[-n]$ are perverse sheaves [BBD], and they naturally underlie mixed Hodge modules (see [Sa1], [Sa2]). These are denoted respectively by $j_!\Q_{h,U}[n]$ and $j_*\D\Q_{h,U}[-n]$ in this paper.
Set $\omega_X(D):=\omega_X\otimes_{\OO_X}\OO_X(D)$ where $\omega_X$ is the (algebraic) dualizing sheaf on $X$. We have the following.
\ms\nin
{\bf Theorem~1.} {\it In the above notation, $\omega_X(D)$ is naturally identified with the first non-zero piece of the Hodge filtration of the underlying filtered right $\DD$-module of $j_*\D\Q_{h,U}[-n]$. Moreover the weight filtration $W$ of $j_*\D\Q_{h,U}[-n]$ induces a canonical weight filtration $W$ on $\omega_X(D)$ such that}
$$\Gr^W_q\omega_X(D)=\mopl_{k+l=n+q+1}\,\omega_{D^{[k,l]}}.$$
\sk
Note that $-q$ coincides with $\dim D^{[k,l]}$ which is equal to $n+1-k-l$. (Here the direct image by the morphism $D^{[k,l]}\to X$ is omitted to simplify the notation.) In case $X$ is smooth, the weight filtration $W$ on $\omega_X(D)$ coincides with the one on $\Omega_X^n(\log D)$ in [De2] up to the shift by $n$ (see also [Fn1]).
Note that the first non-zero piece of the Hodge filtration of the underlying filtered right $\DD$-module is globally well-defined even if $X$ is not globally embeddable into a smooth variety (see (1.1) below).
We can show that the weight filtration $W$ on $\omega_X(D)$ coincides with the one constructed in [Kaw3] since our proof implies the uniqueness of the filtration (see Remark~(2.4) below).
\sk
For a projective morphism $f:X\to Y$, we have the weight spectral sequence of mixed Hodge modules on $Y$
$$E_1^{-q,i+q}=H^if_*\Gr^W_q(j_*\D\Q_{h,U}[-n])\Longrightarrow H^if_*(j_*\D\Q_{h,U}[-n]),$$
degenerating at $E_2$ (see [Sa2]). The restriction of the cohomological direct image of mixed Hodge modules to the first non-zero piece of the Hodge filtration is given by the sheaf-theoretic direct image (see (1.1.3) below). Combining these with Theorem~1 and using the transformation between left and right $\DD_Y$-modules as in (1.2.7) below, we then get the following corollary (which is a slight generalization of the formula (5.4) in [Kaw3]):
\ms\nin
{\bf Corollary~1.} {\it With the notation of Theorem~$1$, let $f:X\to Y$ be a projective morphism of complex algebraic varieties with $Y$ smooth. There is the weight spectral sequence
$$_FE_1^{-q,i+q}=\bigoplus_{k+l=n+q+1}R^if_*\omega_{D^{[k,l]}/Y}\Longrightarrow R^if_*\omega_{X/Y}(D),$$
degenerating at $E_2$, and its $E_1$-differential $d_1$ splits so that the $_FE_2^{-q,i+q}$ are direct factors of $_FE_1^{-q,i+q}$.}
\ms
The last assertion follows from semisimplicity of polarizable Hodge modules. (The latter property comes from that of polarizable variations of pure Hodge structures since polarizable Hodge modules are uniquely determined by their generic variations of Hodge structure.) The above spectral sequence coincides with the one in [Kaw3] obtained under some additional hypotheses as in Corollary~2 below. However, its construction using the residue complex in {\lc} seems rather complicated, and moreover the topological dualizing complex $\D\Q_U$ is not used there so that the definition of the rational structure is quite technical. Here the theories of perverse sheaves [BBD] (see e.g., p.~80--81) and mixed Hodge modules can be used to simplify some arguments (see Remarks~(2.4) below).
\sk
By Corollary~1, the proof of the following semi-positivity theorem for $R^if_*\omega_{X/Y}(D)$ in [FF] can be reduced to the case $X$ is smooth and $D$ is empty, where the assertion has been studied by many people (see also [Fn1], [Fn2], [Fn3], [Ft], [Kaw1], [Kaw2], [Kaw4], [Ko1], [Ko2], [Ko3], [Ko4], [Ko5], [Na]).
\ms\nin
{\bf Corollary~2} ([FF], [Kaw3]). {\it With the notation of Corollary~$1$, assume $Y$ is complete, every $D^{[k,l]}$ is dominant over $Y$, and there is a normal crossing divisor on $Y$ such that the restrictions of the $R^{r-i}(fj)_!\Q_U$ to its complement $Y'$ are locally constant and moreover their local monodromies around it are unipotent where $r=\dim X-\dim Y$. Then $R^if_*\omega_{X/Y}(D)$ is locally free and semi-positive.}
\ms
Here a locally free sheaf $\LL$ on a smooth complete variety $Y$ is called semi-positive if for any morphism $g$ from a smooth complete curve to $Y$, any quotient line bundle (i.e., invertible sheaf) of $g^*\LL$ has non-negative degree (see for instance [Ft], [Kaw1]).
In Corollary~2, it is not necessary to assume that the divisor on $Y$ has {\it simple} normal crossings since it is not needed in Theorems~2 and 3 below.

\sk
Using the above weight spectral sequence together with the generalization of Poincar\'e duality to mixed Hodge modules, we also get a proof of the following assertion in [FF] (see (3.3) below):
\ms\nin
{\bf Theorem~2} ([FF]). {\it With the notation of Corollary~$2$, assume the hypotheses there except for the completeness of $Y$ and the unipotence of the local monodromies. Then $R^if_*\omega_{X/Y}(D)$ coincides with the canonical Deligne extension of the first nonzero piece of the Hodge filtration of the dual variation of mixed Hodge structure of $R^{r-i}(fj)_!\Q_U|_{Y'}$ if $R^if_*\omega_{X/Y}(D)$ does not vanish.}
\ms
Here the canonical extension of the Hodge filtration is defined by taking the intersection of the open direct image with the canonical Deligne extension of the ambient locally free sheaf with an integrable logarithmic connection where the eigenvalues of the residues are contained in $(-1,0]$ (see [De1]). This is closely related with the compatibility of the Hodge filtration $F$ with the filtration $V$ of Kashiwara [Kas1] and Malgrange [Ma] (see [Sa1, Prop.~3.2.2(i)] and also (3.2) below).
We also explain a direct proof of Theorem~2 using the theory of du Bois singularities without using Theorem~1 (see (3.5) below).
Note that Theorem~2 is well-known in case $X$ is smooth and $D=\emptyset$ (see [Ko2], [Na], and also [Sa4, 3.1.1--2]). The proof of Theorem~2 given in [FF] uses a generalization of the argument in [Na] to the mixed case. \sk
Concerning [Kaw3], it is not shown there that the associated variation of mixed Hodge structure is an admissible one, nor that $R^if_*\omega_{X/Y}(D)$ coincides with the canonical extension as is stated in a remark after Thm.~1.1 in {\lc} In fact, for the proof of the semi-positivity, it is sufficient to show this coincidence only for the $E_1$-term of the above spectral sequence which is related to variations of pure Hodge structure (since the semi-positivity is stable by extensions and direct factors in the category of locally free sheaves on $Y$). Here it is enough to assume the unipotence of the local monodromies only for the $E_2$-term, and not for the $E_1$-term, as long as the argument as in the proof of [Kaw1, Thm.~5] is used (although Corollary~2 does not hold without assuming the unipotence of the local monodromies even in the case $X$ is smooth and $D=\emptyset$; see [FF]). Hence the hypothesis in [Kaw3, Thm.~1.1] can be weakened as in our Corollary~2 (unless a geometric argument as in the proof of [Fn2, Thm.~5.4] is used there).
\sk
The proof of Corollary~2 is reduced in our paper by using Theorem~2 to the following semi-positivity theorem for the Hodge filtration of admissible variations of mixed Hodge structure in the sense of [Kas2], [SZ].
\ms\nin
{\bf Theorem~3} ([FF]). {\it Let $Y$ be a smooth complete complex algebraic variety, and $E$ be a normal crossing divisor on it. Let $(M,F)$ be a filtered vector bundle on $Y$. Assume its restriction to the complement of $E$ underlies an admissible variation of mixed $\R$-Hodge structure whose local monodromies around $E$ are unipotent, and $(M,F)$ is its Deligne extension. Let $p_0$ be the maximal integer with $F^pM\ne 0$. Then $F^{p_0}M$ is semi-positive.}
\ms
This is a generalization of a theorem of S.~Zucker [Zu] in the curve case which uses a formula for the curvature due to P.~ Griffiths [Gr].
The proof of Theorem~3 in the mixed case is easily reduced to the pure case by using the admissibility of the variation (see (4.5) below), and moreover the latter case is enough for the proof of Corollary~2 as is explained above.
To simplify the argument, it is very important to proceed by induction on $\dim Y$, and not on the dimension of the stratum strictly containing the image of the curve.
\sk
Essentially the same idea was already used in the proof of [Kaw1, Thm.~5].
Note, however, that the admissibility of the variation of mixed Hodge structure is essential also for the inductive argument there (contrary to the case of [Kaw3, Thm.~1.1]), although it was obtained as a corollary of the multi-variable ${\rm SL}_2$-orbit theorem in [CKS] which was written after [Kaw1]; see also Remark~(4.6)(i) below. (Some remarks related to papers of Kawamata may be found in Remarks~(2.4), (4.4), and (4.6).)
\sk
Note finally that the main theorems in this paper hold also in the analytic case (see Remark~(4.7) below).
\ms\nin
{\bf Acknowledgments.}
The first author was partially supported by the Grant-in-Aid for Young Scientists (A) \#24684002 from JSPS.
The third author is partially supported by Kakenhi 24540039.
The first author would like to thank Professors Takeshi Abe, Hiraku Kawanoue, Kenji Matsuki, and Shigefumi Mori for discussions. He also thanks Professors Valery Alexeev and Christopher Hacon for introducing him this problem. He thanks Professor Gregory James Pearlstein and Professor Fouad El Zein for answering his questions and giving him some useful comments.
The first and second authors would like to thank Professors Kazuya Kato, Chikara Nakayama, and Sampei Usui for fruitful discussions and comments.
The third author thanks Professor Steven Zucker for useful comments on earlier versions of this paper. Finally the authors thank the referee for useful comments to improve the paper.
\sk
This paper is organized as follows.
In Section~1 we review some basic facts in the theory of mixed Hodge modules.
In Section~2 we prove Theorem~1 by calculating the extension classes.
In Section~3 we give a proof of Theorem~2 using Theorem~1 together with another proof using the theory of du Bois singularities.
In Section~4 we give a simple proof of Theorem~3 by using [CK], [CKS], [Sch], [Zu] (without using mixed Hodge modules).
\bs\bs
\centerline{\bf 1. Dualizing sheaves and mixed Hodge modules}
\bs\nin
In this section we review some basic facts in the theory of mixed Hodge modules. Here we use algebraic coherent sheaves and algebraic $\DD$-modules.
\ms\nin
{\bf 1.1.~First non-zero piece of the Hodge filtration.}
Let $(M,F)$ be the underlying filtered $\DD$-module of a mixed Hodge module $\M$ on a complex algebraic variety $X$. This is represented by a system of filtered right $\DD$-modules $(M_{U\into V},F)$ for closed immersions $U\into V$ where $U$ is an open subvariety of $X$ and $V$ is a smooth variety, and they satisfy a certain compatibility condition (see [Sa1, 2.1.20] and [Sa2, 2.1]). If $X$ is globally embeddable into a smooth variety, e.g., if $X$ is projective, then we may assume $U=X$. Set
$$p(\M)=p(M,F):=\min\{p\in\Z\mid F_pM_{U\into V}\ne 0\,\,\,\h{for some}\,\,\,U\into V\}.
\leqno(1.1.1)$$
Then $F_{p(\M)}M_{U\into V}$ is independent of $V$, and depends only on $(M,F)$ and $U$. Indeed, the direct image of a filtered right $\DD_V$-module $(M_V,F)$ under a closed immersion $i:V\into V'$ of smooth varieties of codimension $r$ is locally given by
$$i_*(M_V,F)=(M_V[\xi],F)\q\h{with}\q F_p(M_V[\xi])=\msum_{\nu\in\N^r}\,F_{p-|\nu|}M_V\,\xi^{\nu},
\leqno(1.1.2)$$
where $M_V[\xi]:=M_V[\xi_1,\dots,\xi_r]$ with $\xi:=\partial/\partial x_i$ if $(x_1,\dots,x_m)$ is a local coordinate system of $V'$ with $V=\mcap_{i\le r}\{x_i=0\}$. (More precisely, the direct image of right $\DD_V$-modules by $i:V\into V'$ is defined by the tensor with $\OO_V\otimes_{\OO_{V'}}\DD_{V'}$.) Using the compatibility condition, we can verify that they give a globally well-defined $\OO_X$-module which will be denoted by
$$F_{p(\M)}(\M)\q\h{or}\q F_{p(\M)}(M,F).$$
Here we have to use right $\DD$-modules to get a globally well-defined $\OO$-module (and we need [Sa1, Lemma~3.2.6] to show that $F_{p(\M)}M_{U\into V}$ is an $\OO_U$-module).
\sk
Let $f:X\to Y$ be a projective morphism of complex algebraic varieties.
By the definition of the cohomological direct image of filtered $\DD$-modules $\HH^jf_*$, we get
$$R^jf_*F_{p(\M)}M=F_{p(\M)}\HH^jf_*(M,F)=F_{p(\M)}(H^jf_*\M).
\leqno(1.1.3)$$
See also [Sa4]. Here we have $p(H^jf_*\M)\ge p(\M)$, and the second and the last terms of (1.1.3) are defined to be zero in case $p(H^jf_*\M)>p(\M)$. For the proof of (1.1.3) we may assume $X,Y$ smooth since the assertion is local on $Y$ and $f$ is projective. Using the factorization of $f$ by the composition of the closed immersion $X\into X\times Y$ and the projection $X\times Y\to Y$, the assertion is reduced to the projection case. Here the direct image is defined by using the relative de Rham complex $\DR_{X\times Y/Y}(M,F)$ with $j$-th component given by
$$\h{$\bigwedge $}^{-j}\Theta_X\otimes_{\OO_X}(M,F[-j])\q(j\le 0),
\leqno(1.1.4)$$
where $(F[-j])_p:=F_{p+j}$ (and $pr^{-1}_1$ is omitted before $\Theta_X$ and $\OO_X$). So the first isomorphism of (1.1.3) follows. For the last isomorphism of (1.1.3) we use the fact that the cohomological direct image functor of filtered $\DD$-modules $\HH^jf_*$ is compatible with the cohomological direct image functor of mixed Hodge modules $H^jf_*$. (This is not completely trivial since the direct image of mixed Hodge modules $f_*$ is defined by using Beilinson's resolution.)
\sk
Note that (1.1.3) is a special case of the isomorphisms
$$R^jf_*\Gr^F_p\DR(M)=\HH^j\Gr^F_p\DR(f_*(M,F))=\HH^j\Gr^F_p\DR(f_*(\M)).
\leqno(1.1.5)$$
which follows from the isomorphism $f_*\DR(M,F)=\DR(f_*(M,F))$ (see [Sa1, 2.3.7]).
\ms\nin
{\bf 1.2.~Dualizing sheaves.}
Let $Y$ be a smooth complex variety of dimension $m$. Then $(\omega_Y,F)$ is the underlying filtered right $\DD_Y$-module of the mixed Hodge module
$$\D\Q_{h,Y}[-m]=(\Q_{h,Y}[m])(m),
\leqno(1.2.1)$$
which is pure with weight $-m$. Here $F$ on $\omega_Y$ is defined by
$$\Gr^F_p\omega_Y=0\q(p\ne 0),
\leqno(1.2.2)$$
so that
$$p(\D\Q_{h,Y}[-m])=0,\q F_0(\D\Q_{h,Y}[-m])=\omega_Y.
\leqno(1.2.3)$$
\sk
On the other hand, $(\Omega^m_Y,F)$ is the underlying filtered right $\DD_Y$-module of $\Q_{h,Y}[m]$ which is pure with weight $m$. Here $F$ on $\Omega^m_Y$ is defined by
$$\Gr^F_p\Omega^m_Y=0\q(p\ne -m),
\leqno(1.2.4)$$
so that
$$p(\Q_{h,Y}[m])=-m,\q F_{-m}(\Q_{h,Y}[m])=\Omega^m_Y,
\leqno(1.2.5)$$
and
$$\q(\omega_Y,F)=(\Omega^m_Y,F)(m),
\leqno(1.2.6)$$
This is compatible with (1.2.1).
\sk
Note that $(\Omega^m_Y,F)$ is used for the transformation between filtered left and right $\DD$-modules in [Sa1], [Sa2]; it is defined by associating the following filtered right $\DD_Y$-module to a filtered left $\DD_Y$-modules $(M,F)$:
$$(\Omega^m_Y,F)\otimes_{\OO_Y}(M,F).
\leqno(1.2.7)$$
This transformation is also expressed by choosing local coordinates $y_1,\dots,y_m$ of $Y$ and using the anti-involution $^*$ of $\DD_Y$ defined by
$$(PQ)^*=Q^*P^*,\q (\partial/\partial y_i)^*=-\partial/\partial y_i,\q g^*=g\,\,\,\h{for}\,\,g\in\OO_Y.
\leqno(1.2.8)$$
Here $\Omega_Y^m$ is trivialized by using $dy_1\wedge\cdots\wedge dy_m$.
This expression follows from the definition of the action of $\DD_Y$ on the right $\DD$-module $\Omega_Y^m$ as is well-known.
\ms\nin
{\bf 1.3.~Extension groups.}
Let $X$ be a smooth complex variety, and $Y$ be a smooth closed subvariety of $X$ with $i:Y\into X$ the natural inclusion. Set $c={\rm codim}_YX$. Then
$$\Ext^1_X(\omega_Y,\omega_X)=\Ext^1_Y(\omega_Y,i^!\omega_X)=\Ext^{1-c}_Y(\omega_Y,\omega_Y).
\leqno(1.3.1)$$
Let $\D$ denote the dualizing functor in the bounded derived category of coherent sheaves. This commutes with the direct image $i_*$ by Grothendieck duality. We have
$$\D(\OO_X[\dim X])=\omega_X,\q\D(\OO_Y[\dim Y])=\omega_Y,$$
and the above extension groups are canonically isomorphic to
$$\Ext^{1-c}_X(\OO_X,\OO_Y)=\Ext^{1-c}_Y(i^*\OO_X,\OO_Y)=\Ext^{1-c}_Y(\OO_Y,\OO_Y).
\leqno(1.3.2)$$
\sk
If $c=1$, there are canonical elements $e^*_{X,Y}$ and $e_{X,Y}$ (up to a sign) respectively in the first term of (1.3.1) and (1.3.2). They are dual to each other (up to a sign), and correspond respectively to the short exact sequences
$$0\to\OO_X(-Y)\to\OO_X\to\OO_Y\to 0,
\leqno(1.3.3)$$
$$0\to\omega_X\to\omega_X(Y)\buildrel{\rm res}\over\longrightarrow\omega_Y\to 0,
\leqno(1.3.4)$$
where ``res" is the residue morphism. Here some sign may occur since (1.3.3) corresponds to the distinguished triangle in the derived category of coherent sheaves
$$\OO_Y[\dim Y]\to\OO_X(-Y)[\dim X]\to\OO_X[\dim X]\buildrel{+1\,}\over\longrightarrow,
\leqno(1.3.5)$$
where $+1$ over the last arrow means that the next complex, which is omitted to simplify the notation, is shifted by $+1$.
\sk
We have a similar calculation of extension groups with $\OO_X[n]$ and $\omega_X$ replaced by the constant sheaf $\Q_X[n]$ and its dual $\D\Q_X[-n]$, where $n:=\dim X$.
This can be extended to the case of the mixed Hodge modules. Setting $U:=X\setminus Y$ with $j:U\into X$ the natural inclusion, we get short exact sequences of mixed Hodge modules which are dual to each other:
$$0\to\Q_{h,Y}[n-1]\to j_!\Q_{h,U}[n]\to\Q_{h,X}[n]\to 0,
\leqno(1.3.6)$$
$$0\to\D\Q_{h,X}[-n]\to j_*\D\Q_{h,U}[-n]\to\D\Q_{h,Y}[1-n]\to 0.
\leqno(1.3.7)$$
\sk
Note that we get (1.3.4) by restricting (1.3.7) to the first nonzero piece of the Hodge filtration, i.e., we have in the notation (1.1)
$$p(\D\Q_{h,X}[-n])=p(j_*\D\Q_{h,U}[-n])=p(\D\Q_{h,Y}[1-n])=0,$$
and the short exact sequence (1.3.4) can be identified with
$$0\to F_0(\D\Q_{h,X}[-n])\to F_0(j_*\D\Q_{h,U}[-n])\to F_0(\D\Q_{h,Y}[1-n])\to 0.
\leqno(1.3.8)$$
\sk
This shows that the extension class between the mixed Hodge modules is compatible with the corresponding one between the dualizing sheaves by restricting to $F_0$.

\bs\bs
\centerline{\bf 2. Proof of Theorem~1}
\bs\nin
In this section we prove Theorem~1 by calculating the extension classes. Here we use algebraic coherent sheaves and algebraic $\DD$-modules.
\ms\nin
{\bf 2.1.~Dualizing sheaves of simple normal crossing varieties.}
Let $X$ be a simple normal crossing variety with $X^{[k]}$ as in the introduction. There is a Cech complex
$$\CC_X^{\ssb}:=\bl[\OO_{X^{[1]}}\to\OO_{X^{[2]}}\to\cdots\br]\q\h{with}\q\CC_X^k:=\OO_{X^{[k+1]}},$$
together with a natural quasi-isomorphism
$$\OO_X\simto\CC_X^{\ssb},
\leqno(2.1.1)$$
by choosing an order of the set of irreducible components of $X$.
Set $n=\dim X$. Let $W$ be the increasing filtration on $\CC_X^{\ssb}[n]$ defined by the truncations $\sigma_{\ge -k}$ in [De2] so that
$$\Gr^W_k(\CC_X^{\ssb}[n])=\OO_{X^{[n+1-k]}}[k]\q(k\in[0,n]).
\leqno(2.1.2)$$
Here note that $\dim X^{[n+1-k]}=k$. Then $W$ induces the dual (decreasing) filtration $W$ on$$\omega_X=\D(\OO_X[n])=\D(\CC_X^{\ssb}[n]),$$
satisfying
$$\Gr_W^k\,\omega_X=\omega_{X^{[n+1-k]}}\q(k\in[0,n]).
\leqno(2.1.3)$$
\sk
We have a similar construction with $\OO_X[n]$ and $\omega_X$ replaced by the constant sheaf $\Q_X[n]$ and its dual $\D\Q_X[-n]$, where we get filtrations of perverse sheaves since the shifted Cech complex endowed with the truncations $\sigma_{\ge -k}$ defines a filtration of perverse sheaves (see [BBD]).
This can be extended to the case of the mixed Hodge modules $\Q_{h,X}[n]$ and $\D\Q_{h,X}[-n]$ where we get filtrations of mixed Hodge modules.
\ms\nin
{\bf 2.2.~Dualizing sheaves of reduced simple normal crossing pairs.}
Let $(X,D)$ be a reduced simple normal crossing pair, and $X^{[k]}$, $D^{[k]}$, $D^{[k,l]}$ be as in the introduction. Since $D^{[k]}$ is a simple normal crossing divisor on $X^{[k]}$ for each $k>0$, there is a Cech complex
$$\CC_{X^{[k]},D^{[k]}}^{\ssb}:=\bl[\OO_{X^{[k]}}\to\OO_{D^{[k,1]}}\to\cdots\br]\q\h{with}\q\CC_{X^{[k]},D^{[k]}}^l:=\OO_{D^{[k,l]}},$$
together with a natural quasi-isomorphism
$$\OO_{X^{[k]}}(-D^{[k]})\simto\CC_{X^{[k]},D^{[k]}}^{\ssb}.
\leqno(2.2.1)$$
These give a complex $\CC_{X,D}^{\ssb}$ with
$$\CC_{X,D}^q:=\mopl_{\dim D^{[k,l]}=q}\,\OO_{D^{[k,l]}},$$
where $k>0$ and $l\ge 0$.
This can be constructed by using the orientation sheaf in [De2], and this sheaf can be trivialized by choosing an order of the set of irreducible components of $D^{[k]}$. However, we might get some complicated signs here unless $(X,D)$ globally comes from a pair of simple normal crossing divisors on a smooth variety, since it is unclear whether there are always (partial) orderings of the sets of irreducible components of the $D^{[k]}$ in a compatible way with the inclusions between them. (In particular, it is not clear if we can always get a co-semi-simplicial complex.)
\sk
We have moreover a natural quasi-isomorphism
$$\OO_X(-D)\simto\CC_{X,D}^{\ssb}.
\leqno(2.2.2)$$
We then get a finite increasing filtration $W$ on $\CC_{X,D}^{\ssb}[n]$ with
$$\Gr^W_q(\CC_{X,D}^{\ssb}[n])=\mopl_{\dim D^{[k,l]}=q}\,\OO_{D^{[k,l]}}[q]\q(q\in[0,n]).
\leqno(2.2.3)$$
This induces a decreasing filtration $W$ on
$$\omega_X(D)=\D(\OO_X(-D)[n])=\D(\CC_{X,D}^{\ssb}[n])$$
such that
$$\Gr_W^q\,\omega_X(D)=\mopl_{\dim D^{[k,l]}=q}\,\omega_{D^{[k,l]}}\q(q\in[0,n]).
\leqno(2.2.4)$$
\sk
The above argument can be carried out for the corresponding mixed Hodge modules by replacing $\OO_X(-D)[n]$ and $\omega_X(D)$ with $j_!\Q_{h,U}[n]$ and $j_*\D\Q_{h,U}[-n]$ as in the last remarks of (1.3) and (2.1). Here the shifted Cech complex endowed with the truncations $\sigma_{\ge -k}$ becomes a mixed Hodge module with the weight filtration $W$. We then get decompositions of pure Hodge modules on $X$
$$\Gr^W_q(j_!\Q_{h,U}[n])=\mopl_{\dim D^{[k,l]}=q}\,\Q_{h,D^{[k,l]}}[q]\q(q\in[0,n]),
\leqno(2.2.5)$$
$$\Gr^W_{-q}(j_*\D\Q_{h,U}[-n])=\mopl_{\dim D^{[k,l]}=q}\,\D\Q_{h,D^{[k,l]}}[-q]\q(q\in[0,n]),
\leqno(2.2.6)$$
which are dual of each other. Here $W^q=W_{-q}$ and $\Gr_W^q=\Gr^W_{-q}$ as usual.
\ms
We can now prove the following theorem which implies Theorem~1 in the introduction.
\ms\nin
{\bf 2.3.~Theorem.} {\it In the notation of the introduction and $(1.1)$, we have $p(j_*\D\Q_{h,U}[-n])=0$ and there is a canonical filtered isomorphism of $\OO_X$-modules}
$$F_0(j_*\D\Q_{h,U}[-n],W)=(\omega_X(D),W).
\leqno(2.3.1)$$
\ms\nin
{\it Proof.} The first assertion $p(j_*\D\Q_{h,U}[-n])=0$ follows from (2.2.6). Set
$$(A,W)=(\omega_X(D),W),\q (B,W)=F_0(j_*\D\Q_{h,U}[-n],W).$$
After passing to the graded pieces $\Gr_W^q$, we have the following canonical isomorphisms by (2.2.4), (2.2.6) and (1.2.3):
$$\Gr_W^qA=\Gr_W^qB\q(q\in[0,n]).
\leqno(2.3.2)$$
By induction on $q$, we show the following canonical isomorphisms lifting (2.3.2):
$$(W^qA,W)=(W^qB,W)\q(q\in[0,n]).
\leqno(2.3.3)$$
\sk
The assertion holds for $q=n$ by (2.3.2) with $q=n$. Let $q\in[0,n-1]$, and assume (2.3.3) holds for $q+1$. Consider the extension classes
$$e_A^q\in\Ext^1\bl(\Gr_W^qA,\Gr_W^{q+1}A\br),\q e_B^q\in\Ext^1\bl(\Gr_W^qB,\Gr_W^{q+1}B\br),$$
associated with $(W^qA/W^{q+2}A,W)$ and $(W^qB/W^{q+2}B,W)$. These are identified with each other under the isomorphisms (2.3.2), since the short exact sequence (1.3.4) is identified with (1.3.8). (Note that the signs coming from the orientation sheaf in [De2] are the same for $A$ and $B$.)
\sk
Consider the short exact sequences
$$0\to W^{q+2}A\to W^{q+1}A\to\Gr_W^{q+1}A\to 0.$$
We have the canonical inclusion
$$\Ext^1(\Gr_W^qA,W^{q+1}A)\into\Ext^1(\Gr_W^qA,\Gr_W^{q+1}A),
\leqno(2.3.4)$$
(and similarly for $B$), since (1.3.1) for $c\ge 2$ implies
$$\Ext^1(\Gr_W^{q}A,W^{q+2}A)=0.$$
Consider the extension classes
$$e_A^{\prime\,q}\in\Ext^1\bl(\Gr_W^qA,W^{q+1}A\br),\q e_B^{\prime\,q}\in\Ext^1\bl(\Gr_W^qB,W^{q+1}B\br),$$
associated with $(W^qA,W)$ and $(W^qB,W)$. These are identified with each other under the isomorphisms (2.3.2) for $q$ and (2.3.3) for $q+1$ by the injectivity of (2.3.4). This means that we get a commutative diagram
$$\begin{CD}0@>>>W^{q+1}A@>>>W^qA@>>>\Gr_W^qA@>>>0\\ @. @| @| @|\\
0@>>>W^{q+1}B@>>>W^qB@>>>\Gr_W^qB@>>>0\end{CD}
\leqno(2.3.5)$$
Here the middle vertical morphism is unique if we have
$${\rm Hom}(\Gr_W^qA,W^{q+1}A)=0.$$
But (2.2.4) implies
$${\rm Hom}\bl(\Gr^i_WA,\Gr^j_WA\br)=0\q\h{for}\,\,\,i<j.$$
So (2.3.3) holds also for $q$. This finishes the proof of Theorem~1.
\ms\nin
{\bf 2.4.~Remarks.} (i) The above proof of Theorem~1 gives the uniqueness of the filtration $W$ on $\omega_X(D)$. Hence it coincides with the filtration constructed in [Kaw3] in the absolute case by using the double complex of logarithmic de Rham complexes with second differential given by residue morphisms. In {\lc}, however, only the components of the graded pieces of the weight filtration are given without caring about the induced differential between them, and it seems necessary to show that each graded piece of the double complex is isomorphic to $\mopl_Z\,\sigma_{\ge 0}(K_Z^{\ssb}[1])$ up to a shift of double complexes, where $K_Z^{\ssb}$ is the Koszul complexes associated with the identity morphisms on the de Rham complexes of closed strata $Z$ of a fixed dimension, and $\sigma_{\ge 0}$ is the truncation as in [De2] with respect to the differential of the Koszul complex so that $\sigma_{\ge 0}(K_Z^{\ssb}[1])$ is quasi-isomorphic to the de Rham complex of $Z$. (It is not sufficient to count the multiplicities of the de Rham complexes appearing in the graded pieces of the weight filtration without calculating the differential.)
\sk
It seems also nontrivial to construct a natural quasi-isomorphism between $\omega_X(D)$ and the highest degree part with respect to the de Rham differential of the double complex mentioned above. With our notation, it seems necessary to construct a canonical morphism
$$\omega_X(D)\to\omega_{X^{[1]}}(D^{[1]}+Z^{[1]}),
\leqno(2.4.1)$$
where $Z^{[1]}$ denotes the intersection of $X^{[1]}$ with pull-back of the singular locus of $X$. (Here the direct image by $X^{[1]}\to X$ is omitted to simplify the notation as in {\lc}) In fact, this can be obtained by using a canonical morphism
$$\OO_{X^{[1]}}(-D^{[1]}-Z^{[1]})\to\OO_X(-D).
\leqno(2.4.2)$$
This is easily reduced to the case $D=D^{[1]}=\emptyset$. Analytic-locally there are coordinates $x_1,\dots,x_{n+1}$ of an ambient smooth variety such that $X$ is locally defined by
$$g:=\mprod_{1\le i\le r}\,x_i.$$
For the construction of (2.4.2) with $D=D^{[1]}=\emptyset$, we may replace $X^{[1]}$ with one component of it which is defined by $x_1=0$. Then $Z^{[1]}$ is defined by
$$g':=\mprod_{2\le i\le r}\,x_i,$$
and the desired morphism is induced by the inclusion of the ideal generated by $g'$ in the structure sheaf of the ambient smooth variety by dividing it by the ideal generated by $g$.
\ms
(ii) For the argument in [Kaw3] (which proves Corollary~2 by showing Corollary~1 in our paper), it is unnecessary to show that the obtained variation of mixed Hodge structure is an admissible one so that the canonical extension of the Hodge filtration is compatible with the passage to the graded pieces of the weight filtration (see the remark before Theorem~3). However, it is possible to show its coincidence with our variation of mixed Hodge structure by using the functor $\DR^{-1}$ as in the remarks below.
\ms
(iii) Assume $X$ is a simple normal crossing divisor on a smooth variety $Y$, and $D$ is the restriction of a normal crossing divisor $D'$ on $Y$ such that $X\cup D'$ is a normal crossing divisor on $Y$. Set
$$U':=Y\setminus D',\q U'':=Y\setminus(X\cup D'),$$
with inclusions $j':U'\into Y$, $j'':U''\into Y$. We have a short exact sequence of mixed Hodge modules
$$0\to j'_*\D_{U'}[-n-1]\to j''_*\D_{U''}[-n-1]\to j_*\D_U[-n]\to 0.
\leqno(2.4.3)$$
In fact, this is easy in case $D'=\emptyset$ so that $U'=Y$, $U''=Y\setminus X$, and $U=X$. In general, apply the functor $j'_*j'{}^*$ to this short exact sequence for the case $D'=\emptyset$.
\sk
The above short exact sequence corresponds to the short exact sequence of coherent sheaves
$$0\to\omega_Y(D')\to\omega_Y(X+D')\to\omega_X(D)\to 0.
\leqno(2.4.4)$$
The corresponding short exact sequence is stated after the definition of the filtered complex in the absolute case [Kaw3, Sect.~3].
We get the underlying short exact sequence of filtered $\DD$-modules of the above short exact sequence of mixed Hodge modules by applying the functor $\DR^{-1}$ to the short exact sequence in {\lc}
\ms
(iv) Kawamata's construction [Kaw3] in the absolute case can be interpreted by using the theories of perverse sheaves and mixed Hodge modules as follows.
\sk
More generally, let $X$ be a complex algebraic variety endowed with a stratification $\{S_k\}$ satisfying the following conditions: each $S_k$ is a locally closed smooth subvariety of $X$ with pure dimension $k$ such that the inclusion $i_k:S_k\into X$ is an affine morphism and
$$\SSS_k\setminus S_k=\mcoprod_{j<k}\,S_j.$$
Let $j_k:X\setminus\SSS_k\into X$ denote the inclusion. We have a decreasing filtration $G$ on $\Q_X$ defined by
$$G^{k+1}\Q_X=j_{k!}\Q_{X\setminus\SSS_k}\q\h{so that}\q \Gr^k_G\Q_X=i_{k!}\Q_{S_k}.$$
Since $S_k$ is smooth and $i_k$ is an affine immersion, $i_{k!}\Q_{S_k}[k]$ is a perverse sheaf on $X$. By the theory of realization functors in [BBD, 3.1], this implies an isomorphism in the derived category of perverse sheaves
$$\Q_X=\bl[\Q_{S_0}\to i_{1!}\Q_{S_1}[1]\to i_{2!}\Q_{S_2}[2]\to\cdots\br],
\leqno(2.4.5)$$
where $i_{k!}\Q_{S_k}[k]$ is put at degree $k$ in the target. By the same argument as in the proof of [Sa5, Prop.~4.1] (where each stratum is assumed affine), this is naturally lifted to an isomorphism in the derived category of mixed Hodge modules
$$\Q_{h,X}=\bl[\Q_{h,S_0}\to i_{1!}\Q_{h,S_1}[1]\to i_{2!}\Q_{h,S_2}[2]\to\cdots\br].
\leqno(2.4.6)$$
Taking the dual, we get an isomorphism
$$\bl[\cdots\to i_{2*}\D\Q_{h,S_2}[-2]\to i_{1*}\D\Q_{h,S_1}[-1]\to\D\Q_{h,S_0}\br]=\D\Q_{h,X}.
\leqno(2.4.7)$$
This corresponds to Kawamata'a complex [Kaw3] in the absolute case with $D=\emptyset$ (in our notation) where the stratification is associated with the normal crossing variety $X$. In general, it is enough to apply the functor $j_*j^*$ to the above complex where $j:X\setminus D\into X$ is the inclusion, since we get the underlying complex of bi-filtered $\DD$-modules of the obtained complex of mixed Hodge modules by applying the functor $\DR^{-1}$ to the bi-filtered complex in [Kaw3] in the absolute case. In fact, we get the underlying bi-filtered $\DD_X$-module of $j_*\Q_{X\setminus D}[\dim X]$ by applying the functor $\DR^{-1}$ to the bi-filtered de Rham complex
$$(\Omega_X^{\ssb}(\log D)[\dim X];F,W),$$
if $D$ is a normal crossing divisor on a smooth variety $X$ with $j:X\setminus D\into X$ as above (see for instance [Sa5, 1.4.5]). Note that if $X$ is smooth, then we have by (1.2.1)
$$\D\Q_X[-\dim X]=(\Q_X[\dim X])(\dim X).$$
We can show the compatibility with the differentials of the complexes by comparing (1.3.4) and (1.3.7) at sufficiently general points of the strata.

\bs\bs
\centerline{\bf 3. Proof of Theorem~2}
\bs\nin
In this section we give a proof of Theorem~2 using Theorem~1 together with another proof using the theory of du Bois singularities. Here we mainly use analytic $\DD$-modules and analytic sheaves by using GAGA (except possibly in (3.4) where algebraic sheaves can also be used).
\ms\nin
{\bf 3.1.~Filtration $V$ of Kashiwara and Malgrange.}
Let $M$ be a regular holonomic left $\DD_Y$-module on a complex manifold $Y$.
Let $E$ be a smooth hypersurface in $Y$.
Let $y$ (resp.\ $\xi_y$) be a locally defined function (resp.\ vector field) such that we have locally
$$D=\{y=0\},\q\langle\xi_y,dy\rangle=1.$$
Let $V_0\DD_Y\subset \DD_Y$ be the subring generated by $\OO_Y$ and the vector fields preserving the ideal of $E\subset Y$. Let $\Sigma$ be a subset of $\C$ containing $0$ and such that the composition of the inclusion $\Sigma\into\C$ with the exponential map: $\C\ni\alpha\mapsto\exp(2\pi i\alpha)\in\C^*$ is bijective.
There is a unique separated exhaustive filtration $V$ of Kashiwara [Kas1] and Malgrange [Ma] along $E$ on $M$ indexed by $\Z$ and satisfying the following conditions:
\sk
\begin{itemize}
\item[(3.1.1)]
$\xi_y(V_kM)\subset V_{k+1}M,\,\,\,y(V_kM)\subset V_{k-1}M\,\,\,(\forall\,k)$,
\sk
\item[(3.1.2)]
$y(V_kM)=V_{k-1}M\q(k<0)$,
\sk
\item[(3.1.3)]
$V_kM$ are coherent $V_0\DD_Y$-modules,
\sk
\item[(3.1.4)]
There is a minimal polynomial for the action of $\xi_yy\in V_0\DD_Y$ on $\Gr^V_kM$ such that its roots are contained in $-k+\Sigma$ for any $k$.
\end{itemize}
\sk\nin
These conditions are independent of the choice of $y,\xi_y$, and $V$ exists globally (although there is no canonical $\DD_E$-module structure on $\Gr^k_VM$ without choosing the function $y$).
We say that $M$ is quasi-unipotent along $E$ if the roots of the minimal polynomial are rational numbers. (This condition is satisfied for mixed Hodge modules.) In this case $V$ can be indexed by $\Q$ by replacing $k$ in conditions (3.1.1--3) with rational numbers $\alpha$, and condition (3.1.4) with
\sk
\begin{itemize}
\item[(3.1.5)]
The action of $\xi_yy+\alpha$ on $\Gr^V_{\alpha}M$ is nilpotent.
\end{itemize}
\sk
Note that $(\xi_yy+\alpha)^*=-y\xi_y+\alpha$ in case $\xi_y=\partial/\partial y$ (by choosing local coordinates containing $y$), where $^*$ is the anti-involution as in (1.2.8).
\sk
Assume $M$ is a regular holonomic $\DD_Y$-module of Deligne type, i.e., there is a simple normal crossing divisor $E$ on $Y$ such that $M|_{Y\setminus E}$ is a locally free sheaf with an integrable connection and $M$ is its meromorphic extension. Assume $M$ is quasi-unipotent along any irreducible component $E_l$ of $E$. Let $V^{(l)}$ be the $V$ filtration along $E_l$ indexed by $\Q$. Then the Deligne canonical extension with eigenvalues of the residues contained in $(\alpha-1,\alpha]$ is given by
$$\mcap_l V^{(l)}_{<-\alpha}M.
\leqno(3.1.6)$$
\ms\nin
{\bf 3.2~Compatibility with the Hodge filtration.} Let $(M,F)$ be the underlying filtered left $\DD_Y$-module of a mixed Hodge module on a smooth complex algebraic variety $Y$. Let $E$ be a smooth divisor on $Y$. In the notation of (3.1) the following conditions are satisfied (see [Sa1, 3.2.1]:
$$y(F_pV_{\alpha}M)=F_pV_{\alpha-1}M\q(\alpha<0),
\leqno(3.2.1)$$
$$\xi_y(F_p\Gr^V_{\alpha}M)=F_{p+1}\Gr^V_{\alpha+1}M\q(\alpha>-1).
\leqno(3.2.2)$$
\sk
Let $j:Y\setminus E\into Y$ denote the inclusion.
By Prop.~3.2.2 in {\lc}, condition (3.2.1) is equivalent to the following:
$$F_pV_{<0}M=V_{<0}M\cap j_*j^*F_pM.
\leqno(3.2.3)$$
In case $E$ is a divisor with normal crossings, (3.2.3) implies in the notation of (3.1.6)
$$F_p\bl(\mcap_l V^{(l)}_{<0}M\br)=\mcap_l V^{(l)}_{<0}M\cap j_*j^*F_pM,
\leqno(3.2.4)$$
where the right-hand side is the canonical extension of the Hodge filtration $F$ by the assertion concerning (3.1.6). Indeed, (3.2.4) follows from (3.2.3) by induction on the number of local irreducible components. Here we have to show the inclusion $\supset$. We apply (3.2.3) to a local irreducible component $E_l$ of the divisor $E$ and the inductive hypothesis is used outside $E_l$.
\sk
If $M$ has no-nontrivial quotient $\DD_Y$-submodule supported on $E$, or equivalently, we have the surjectivity of
$$\xi_y:\Gr^V_{-1}M\to\Gr^V_0M,
\leqno(3.2.5)$$
(see {\lc}, Prop.~3.1.8), then (3.2.2) holds also for $\alpha=-1$, i.e.,
$$\xi_y(F_p\Gr^V_{-1}M)=F_{p+1}\Gr^V_0M.
\leqno(3.2.6)$$
Indeed (3.2.5) underlies a morphism of mixed Hodge modules denoted by Var in {\lc}, 5.1.3.4, and the latter is strictly compatible with the Hodge filtration $F$ by Prop.~5.1.14 in {\lc}
\sk
As a conclusion, if (3.2.5) is surjective, then we have for $p=p(M,F)$
$$F_{p(M,F)}M=V_{<0}M\cap j_*j^*F_{p(M,F)}M.
\leqno(3.2.7)$$
\ms\nin
{\bf 3.3.~Proof of Theorem~2} (using Theorem~1).
By duality of mixed Hodge modules (see [Sa2, 4.3.5]), there is a canonical isomorphism
$$H^if_*(j_*\D\Q_{h,U}[-n])=\D(H^{-i}f_*(j_!\Q_{h,U}[n])).
\leqno(3.3.1)$$
Here the right-hand side is the dual as a mixed Hodge module, and its restriction over $Y'$ is identified (up to certain shifts of filtrations) with the dual of the variation of mixed Hodge structure associated with
$$H^{r-i}(fj)_!\Q_U|_{Y'}.$$
\sk
Let $(M;F,W)$ be the bi-filtered left $\DD_Y$-modules underlying
$$H^if_*(j_*\D\Q_{h,U}[-n]).$$
Here the Hodge filtration $F$ is indexed like a right $\DD$-module without taking the shift of filtration in the transformation between left and right $\DD$-modules, i.e., the tensor with $(\omega_Y,F)$ is used instead of $(\Omega_X^m,F)$ in (1.2.7). So we have in the notation of (1.1.1)
$$p(M,F)=0.
\leqno(3.3.2)$$
\sk
Let $V^{(l)}$ denote the filtration of Kashiwara and Malgrange on $M$ along an irreducible component $E_l$ of $E$.
We first show
$$F_0\Gr^{V^{(l)}}_{\alpha}M=0\q(\alpha\ge 0),\q\h{i.e.,}\q F_0M=F_0V^{(l)}_{<0}M.
\leqno(3.3.3)$$
Indeed, (3.3.3) for $\alpha>0$ follows from (3.2.2) together with (3.3.2).
For $\alpha=0$, we have
$$F_0\Gr^{V^{(l)}}_0\Gr^W_kM=0\q(\alpha\ge 0).
\leqno(3.3.4)$$
This follows from [Sa4, 2.6.1] by using the weight spectral sequence in the introduction, and reducing to the case $X$ is smooth and $D=\emptyset$, since every stratum is dominant over $Y$. Here $F,W,V^{(l)}$ are compatible three filtrations on $M$ (see [Sa2, 2.2.1]). So (3.3.3) follows.
(Note that the surjectivity of (3.2.5) does not necessarily hold on $Y$ although it holds on $X$ by the assumption that any stratum is dominant over $Y$; see (3.5) below.)
\sk
Let $j:Y\setminus E\into Y$ denote the inclusion. We get by (3.3.3) and (3.2.4)
$$F_0M=F_0\bl(\mcap_lV^{(l)}_{<0}M\br)=\mcap_lV^{(l)}_{<0}M\cap j_{*}j^*F_0M.
\leqno(3.3.5)$$
So the assertion follows.
\ms\nin
{\bf 3.4.~Du Bois singularities.} Let $X$ be a reduced complex algebraic variety. Let $(\widetilde{\Omega}_X^{\ssb},F)$ denote the filtered du Bois complex (see [dB]). It is isomorphic to $\DR(\Q_{h,X})$ which is obtained by applying the de Rham functor $\DR$ to the underlying complex of filtered $\DD$-modules of $\Q_{h,X}$. (Here we assume $X$ is embeddable into a smooth variety for simplicity since the problem is local.) In particular, we get
$$\Gr^F_0\DR(\Q_{h,X})=\Gr_F^0\widetilde{\Omega}_X^{\ssb},
\leqno(3.4.1)$$
where $F^p=F_{-p}$. By the definition of $\D$ for filtered differential complexes in [Sa1,2.4.11], $\DR$ is compatible with $\D$, and we get
$$\D\Gr^F_0\DR(\Q_{h,X})=\Gr^F_0\DR(\D\Q_{h,X})=F_0\DR(\D\Q_{h,X}),
\leqno(3.4.2)$$
where the last isomorphism follows from the definition of $\DR$ (see (1.1.4)) together with
$$p(\D\Q_{h,X})=0.
\leqno(3.4.3)$$
The last equality can be reduced to the smooth case by using a smooth affine stratification as in Remark~(2.4)(iv).
\sk
Recall that $X$ is called du Bois if
$$\Gr_F^0\widetilde{\Omega}_X^{\ssb}=\OO_X.$$
By the above argument, this condition is equivalent to
$$F_0\DR(\D\Q_{h,X})=\K_X,
\leqno(3.4.4)$$
where $\K_X$ denotes the dualizing complex of $X$. If $X$ is Cohen-Macaulay, then we have
$$\K_X=\omega_X[\dim X].$$
\sk
Let $D$ be a locally principal reduced divisor on $X$ with $j:X\setminus D\into X$ the canonical inclusion. Assume $X,D$ are du Bois. Then
$$\OO_X(-D)=\Gr^F_0\DR(j_!\Q_{h,X\setminus D}),\q\K_X(D)=F_0\DR(j_*\D\Q_{h,X\setminus D}).
\leqno(3.4.5)$$
Indeed, we get the first isomorphism by comparing the following distinguished triangle and exact sequence:
$$\aligned\Gr^F_0\DR(j_!\Q_{h,X\setminus D})\to\Gr^F_0\DR(\Q_{h,X})&\to\Gr^F_0\DR(\Q_{h,D})\buildrel{+1}\over\longrightarrow,\\0\to\OO_X(-D)\to\OO_X&\to\OO_D\to 0.\endaligned
\leqno(3.4.6)$$
The second isomorphism then follows by using duality together with (3.4.3). (See also [Kov].)
\ms\nin
{\bf 3.5.~Direct proof of Theorem~2.} Let $(X,D)$ be a reduced simple normal crossing pair.
It is well-known that $X$ is du Bois by using Remark~(2.4)(iii) together with the weight filtration. (See also [Kov].)
As for $D$, there are local coordinates $x_1,\dots,x_{n+1}$ in some ambient smooth variety such that we have analytic-locally
$$D=\bl\{\mprod_{1\le i\le r}\,x_i=0\br\}\cap\bl\{\mprod_{r<i\le m}\,x_i=0\br\}\subset\Delta^r\times\Delta^{n+1-r},$$
since $D$ is locally the intersection of two divisors having no common irreducible components and whose union is a normal crossing divisor. So $D$ is analytic-locally a product of two normal crossing divisors, and is du Bois.
\sk
Let $g$ be a local defining equation of an irreducible component of $E\subset Y$. Set $g':=f^*g$. Consider the graph embeddings
$$i_{g'}:X\into X\times\C,\q i_g:Y\into Y\times\C.$$
They are compatible with
$$f:X\to Y,\q f':=f\times id:X\times\C\to Y\times\C.$$
\sk
Let $(M',F)$ be the underlying filtered $\DD$-module of $i_{g*}j_*\D_{h,U}[-n]$. By (3.4.5) and (1.1.2) we have
$$F_0M'=\omega_X(D).
\leqno(3.5.1)$$
By the same argument as in (3.3) and using (1.1.3), it is enough to show
$$F_0\Gr^V_{\alpha}H^if'_*M'=0\q\h{for}\,\,\,\alpha\ge 0.
\leqno(3.5.2)$$
Here we may assume that $X$ is a subvariety of a smooth variety since $f$ is projective and the assertion is local on $Y$.
We have the commutativity of $F_0\Gr^V_\alpha$ with the direct image $H^if'$, since $f'_*(M';F,V)$ is bi-strict (see [Sa1, 3.3.17] and also Remark~(3.6) below). So it is enough to show
$$F_0\Gr^V_{\alpha}M'=0\q\h{for}\,\,\,\alpha\ge 0.
\leqno(3.5.3)$$
Then the assertion is local on $X$. It is reduced to the case $X$ is smooth (but $D\ne\emptyset$) by using the surjective morphisms of the short exact sequences (2.4.3--4) in Remark~(2.4)(iii). It is further reduced to the case $X$ is smooth and $D=\emptyset$ by using the weight filtration $W$. So the assertion follows from (3.2.7) (see also [Sa4]) where the assumption on the surjectivity of (3.2.5) is satisfied since it is assumed that every stratum is dominant over $Y$.
This finishes the second proof of Theorem~2.
\ms\nin
{\bf 3.6.~Remark.} We say that a bi-filtered complex $(K;F,G)$ is bi-strict (see [Sa1, 1.2.2]) if we have the injectivity of the morphisms
$$H^iF_pK\into H^iK,\q H^iG_qK\into H^iK,\q H^iF_pG_qK\into H^iK,$$
together with
$$H^iF_pK\cap H^iG_qK=H^iF_pG_qK\q\h{in}\,\,\,H^iK.$$
By \lc, Prop.~1.2.13, we then get the commutativity of the three functors
$$H^i,\,\,\Gr^F_p,\,\,\Gr^G_q.$$
In fact, let $E$ be the filtration on $K^i$ defined by
$$E_{-1}K^i=0,\,\,\,E_0K^i={\rm Im}\,d_{i-1},\,\,\,E_{1}K^i={\rm Ker}\,d_i,\,\,\,E_2K^i=K^i.$$
By \lc, Prop.~1.2.3(iii), the bi-strictness of $(K;F,G)$ is equivalent to the condition that $E,F,G$ on each $K^i$ are compatible three filtrations in the sense of \lc, 1.1.13.
\sk
If $G_qK=0$ ($q\ll 0$) and $G$ is exhaustive, then the bi-strictness of $(K;F,G)$ is equivalent to the condition that $(K,G)$ and the $\Gr^G_q(K,F)$ are strict (see \lc, Cor.~1.2.10). However, this does not apply to the case of the filtration $V$, and we have to use the completion as in \lc, Lemma~3.3.2.
Note that the completion by $V$ is essentially the $y$-adic completion by condition (3.2.1) above.
\bs\bs
\centerline{\bf 4. Proof of Theorem~3}
\bs\nin
In this section we give a simple proof of Theorem~3 by using [CK], [CKS], [Sch], [Zu] (without using mixed Hodge modules). Here we use analytic sheaves.
\ms\nin
{\bf 4.1~Deligne extensions.}
Let $S$ be a polydisk with coordinates $t_1,\dots,t_n$. Set
$$E:=\mcup_{j=1}^n\{t_j=0\}.$$
Let $L$ be a $\C$-local system on $U:=S\setminus E$ with unipotent monodromies. Let $M$ be the Deligne extension [De1] of $\OO_U\otimes_{\C}L$. Let $(u_i)$ be a basis of multivalued sections of $L$. Then it corresponds to a basis $(\uu_i)$ of $M$ such that
$$\uu_i|_U=\exp\bigl(-\msum_j\,(\log t_j)N_j\bigr)u_i,
\leqno(4.1.1)$$
where $N_j=(2\pi i)^{-1}\log T_j$ with $T_j$ the monodromy around $\{t_j=0\}$ (see [De1] and also [Kaw2]).
This implies that the Deligne extensions are stable by the pull-back under ramified coverings of polydisks in the unipotent monodromy case.
\sk
Set $Z:=\mcap_{i\le r}\{t_i=0\}$ where $r\ge 2$. Let $\pi:S'\to S$ be the blow-up of $S$ along $Z$. It has local coordinates $x'_1,\dots,x'_n$ with
$$\pi^*t_i=\begin{cases}x'_ix'_r&\h{if}\,\,\,i<r,\\ x'_i&\h{if}\,\,\,i\ge r.\end{cases}$$
Substituting this to (4.1.1), we see that $\pi^*M$ coincides with the Deligne extension of the pull-back of $L$ (see also the proof of [Kaw2, Prop.~1]).
\ms\nin
{\bf 4.2.~Deligne extensions of filtered vector bundles.}
The above arguments imply that the Deligne extensions of the underlying filtered vector bundles of admissible variations of mixed Hodge structure with unipotent local monodromies are stable by the pull-backs under ramified coverings or blow-ups as above.
In fact, this is reduced to (4.1) since two filtrations of a vector bundle coincide if they coincide on a dense Zariski-open subset. (Note, however, that a vector bundle in this paper means a locally free sheaf of finite type, and a filtered vector bundle means that its graded pieces are locally free.)
\sk
As a corollary, we get the following (which does not seem to be completely trivial; see Remark~(4.4)(i) below).
\ms\nin
{\bf 4.3.~Proposition.} {\it Let $Y$ be a complex manifold with $E$ a normal crossing divisor on it. Let $\M$ be an admissible variation of mixed Hodge structure on $Y\setminus E$ with unipotent local monodromies around $E$, which is admissible with respect to $Y$. Let $(M,F)$ be the Deligne extension of the underlying filtered vector bundle of $\M$ to $Y$. Let $g:(\Delta,0)\to(Y,E)$ be a morphism from an open disk $\Delta$ such that $g^{-1}(E)=\{0\}$. Then $g^*(M,F)$ coincides with the Deligne extension of the pull-back of $\M$ by $g$.}
\ms\nin
{\it Proof.}
We may replace $Y$ with a smooth center blow-up of $Y$ as in (4.1-2) (by factorizing $g$). Repeating this, we may assume that the curve $g(C)$ is nonsingular near $g(0)$, and moreover there are local coordinates $y_1,\dots,y_n$ of $Y$ such that we have locally
$$E=\{y_1=0\},\q g(C)=\mcap_{i=2}^n\{y_i=0\}.$$
So the assertion follows in case $g$ is a closed embedding, and is reduced to (4.2) in case $g$ is a ramified covering over the image. This finishes the proof of Proposition~(4.3).
\ms\nin
{\bf 4.4.~Remarks.} (i) It does not seem that Proposition~(4.3) can be proved trivially without using, for instance, an embedded desingularization of the image of the curve as is explained above, since the image has the Puiseux expansion in case $\dim Y=2$, and it may have many Puiseux pairs.
\sk
We can prove Proposition~(4.3) also by calculating the residue of the pull-back of the logarithmic connection and using the commutativity of the actions of $y_i\partial/\partial y_i$ in the normal crossing case. Here it is also possible to reduce the assertion to the case of a constant local system of rank 1 by taking the graded pieces of an appropriate filtration on the local system after restricting it to $(\Delta^*)^n\subset Y\setminus E$, since the assertion is local on $Y$. Note that canonical Deligne extensions are stable by extensions of locally free sheaves endowed with logarithmic integrable connections; i.e., the middle term of a short exact sequence is a canonical Deligne extension if the other two terms are, since only the eigenvalues of the residues matter.
\sk
(ii) Proposition~(4.3) can be extended to the case of a morphism $g:(X,D)\to(Y,E)$ where $D,E$ are normal crossing divisors on smooth varieties $X,Y$ such that $D\supset g^{-1}(E)$. Indeed, it is sufficient to restrict to curves on $X$ intersecting $D$ at smooth points, since this implies the coincidence outside a subvariety of codimension 2. The assertion also follows from the second proof of Proposition~(4.3) as is explained above.
\ms\nin
{\bf 4.5. Proof of Theorem~3.} In this proof, we assume more generally that the underlying local system of $\M$ is defined over a subfield $A$ of $\R$. (So $A$ can be $\Q$.)
Since semi-positivity is stable by extensions of free sheaves on $Y$, we may assume that the admissible variation $\M$ on the complement of $E$ is pure.
In case $\dim Y=1$, it is reduced to [Zu] by using (4.2) in the ramified covering case.
\sk
In the general case we proceed by induction on $\dim Y$, and not on the dimension of the stratum strictly containing $g(C)$ (where strictly means that the stratum is the minimal one among the strata containing $g(C)$).
Let $g:C\to Y$ be a morphism from a smooth complete curve. If $g(C)$ is not contained in $E$, then the pull-back of the Deligne extension by $g$ coincides with the Deligne extension of the pull-back of the local system by Proposition~(4.3). So the assertion is reduced to [Zu] as is explained above.
\sk
We now consider the case $g(C)\subset E$. Let $D$ be the normalization of an irreducible component of $E$ containing $g(C)$. Let $\rho:D\to Y$ be the natural morphism. Set
$$D':=D\setminus\rho^{-1}({\rm Sing}\,E).$$
This is identified with a locally closed subvariety of $Y$. Take $0\in D'$. Let $y$ be a local defining equation of $D'\subset Y$ at $0$. Let $U$ be a sufficiently small open polydisk around $0\in Y$ associated with some local coordinate system containing $y$. Here we may assume that $U\cap\rho(D)$ coincides with $U\cap D'$. Set $U':=U\setminus D'$. This is a product of a polydisk and a punctured disk. Set $Y':=Y\setminus E$. Let $L$ be the underlying $A$-local system of the variation of Hodge structure on $Y'$. A multivalued horizontal section $u$ of $L_{U'}$ defines a holomorphic local section $\uu$ of $M_U$ as in (4.1.1), i.e.,
$$\uu|_{U'}=\exp(-(\log y)N)u,
\leqno(4.5.1)$$
where $N:=(2\pi i)^{-1}\log T$ with $T$ the local monodromy of $L$ around $D'$. (Note that $(2\pi i)^{-1}$ corresponds to the Tate twist; see [De2].) This induces an injection
$$\psi_yL_{U'}\into M_{U\cap D'},
\leqno(4.5.2)$$
giving an $A$-structure on $M_{U\cap D'}$. Here $\psi_y$ is the nearby cycle functor, and $M_{U\cap D'}$ denotes the restriction of $M$ to $U\cap D'$ as a locally free sheaf (i.e., as a vector bundle). Moreover this is compatible with the induced connection on $M_{U\cap D'}$ (which may depend on the choice of the local coordinate $y$). We have also the induced polarization
$$\psi_yL_{U'}\otimes_{A}\psi_yL_{U'}\to A(-w),
\leqno(4.5.3)$$
by a polarization of the variation of Hodge structure $L\otimes_{A}L\to A(-w)$, where $w$ is the weight of the variation of Hodge structure.
\sk
By (4.5.1) the action of $N$ on the left-hand side of (4.5.2) corresponds to $\theta_y:=-y\partial/\partial y$ on the right-hand side. (Here $\partial/\partial y$ may depend on the coordinate system containing $y$ although it is not clear from the notation.) Moreover the action of $\theta_y$ on $M_{U\cap D'}$ coincides with the residue of the logarithmic connection up to a sign, and is independent of the choice of the coordinates. So it will be denoted also by $N$.
\sk
Let $W$ be the finite increasing filtration on $\psi_yL_{U'}$ and $M_{U\cap D'}$ satisfying
$$\aligned N(W_k\psi_yL_{U'})&\subset W_{k-2}\psi_yL_{U'}\q(k\in\Z),\\ N^k:\Gr^W_{w+k}\psi_yL_{U'}&\simto(\Gr^W_{w-k}\psi_yL_{U'})(-k)\q(k>0),\endaligned
\leqno(4.5.4)$$
and similarly with $\psi_yL_{U'}$ replaced by $M_{U\cap D'}$.
These are compatible with the inclusion (4.5.2). The induced polarization (4.5.3) gives a perfect pairing
$$\Gr^W_{w+k}\psi_yL_{U'}\otimes\Gr^W_{w-k}\psi_yL_{U'}\to A(-w).$$
Combined with the isomorphism in (4.5.4), this induces a polarization on the $N$-primitive part
$$P_N\Gr^W_{w+k}\psi_yL_{U'}:={\rm Ker}\,N^{k+1}\subset\Gr^W_{w+k}\psi_yL_{U'}.$$
We thus get an $A$-structure together with the induced polarization on the $N$-primitive part
$$P_N\Gr^W_{w+k}M_{U\cap D'}:={\rm Ker}\,N^{k+1}\subset\Gr^W_{w+k}M_{U\cap D'}.$$
(Here $P_N\Gr^W_{w+k}\psi_yL_{U'}=P_N\Gr^W_{w+k}M_{U\cap D'}=0$ for $k<0$.)
\sk
We can show that the obtained $A$-structure and polarization on $P_N\Gr^W_{w+k}M_{U\cap D'}$ are independent of the choice of $y$ by the same argument as in [St, 4.24]. (Indeed, this easily follows from (4.5.1) since the action of $N$ on the graded pieces of $W$ vanishes.) So they exist globally on $D'$. By the Lefschetz decomposition
$$\Gr^W_kM_{U\cap D'}=\mopl_{j\ge 0}\,N^j(P_N\Gr^W_{k+2j}M_{U\cap D'}),
\leqno(4.5.5)$$
these give a canonical $A$-structure together with a global polarization on $\Gr^W_k(M,F)_{D'}$ in a compatible way with the induced connection on $\Gr^W_kM_{D'}$. Here the latter is also well-defined since the action of $N$ on $\Gr^W_kM_{D'}$ vanishes.
They form a variation of polarized Hodge structure on $D'$ by [Sch] as is well-known.
Moreover, by the multi-variable ${\rm SL}_2$-orbit theorem, $(M_{D'};F,W)$ locally underlies an admissible variation of mixed Hodge structure with respect to the compactification $D\supset D'$ by taking local coordinates compatible with $E$ (see [CK], [CKS]).
Here $(M,F)_D:=\rho^*(M,F)$ is the Deligne extension of $(M,F)_{D'}$ since this holds by forgetting $F$.
Then $\Gr^W_k(M,F)_D$ is the Deligne extension of $\Gr^W_k(M,F)_{D'}$ where $W$ on $M_D$ is defined by using the Deligne extension (see Remark~(4.6)(i) below).
So we can proceed by induction on $\dim Y$, replacing $Y$ with $D$, and $(M,F)$ with $\Gr^W_k(M,F)_D$.
This finishes the proof of Theorem~3.
\ms\nin
{\bf 4.6.~Remarks.} (i) In the above proof, the key point is the following:
$$\h{$\Gr^p_F\Gr^W_kM_D$ are locally free for any $p,k$.}
\leqno(4.6.1)$$
This is one of the conditions for admissible variation of mixed Hodge structure in the curve case (see [SZ]), and follows from the multi-variable ${\rm SL}_2$-orbit theorem [CKS]. Without this property, it is not clear whether the Deligne extension of the Hodge filtration is compatible with the passage to the graded pieces of the weight filtration $W$. The last property is essential for the inductive argument in the proof of Theorem~3.
This is also the most nontrivial point in the proof of [Kaw1, Thm.~5].
\ms
(ii) Only a sketch of the proof was provided for the case of strata with codimension 1 or 2 in \lc, where the structure of the induction is not very clear: If one uses the induction on $\dim Y$, then the classification by the dimension of the stratum strictly containing the image of the curve $g(C)$ is unnecessary, and it is enough to consider only the two cases: $g(C)\subset E$ or not. If one uses the induction on the dimension of the stratum strictly containing the image of the curve, then we have to treat the case of higher codimensional case where the argument becomes much more complicated. (In fact, this is closely related to the proof of the theorem that admissible variations of mixed Hodge structures are mixed Hodge modules in [Sa2, 3.c], where the compatibility between the induced dualities of the nearby cycles for constructible sheaves and for filtered $\DD$-modules is proved in a more general situation.)
\sk
Another minor problem in the proof of [Kaw1, Thm.~5] is that the definition of the induced polarization on the graded pieces of the weight filtration was not very precise. In fact, we have to use a local coordinate for this, and then prove the independence of the choice of the coordinate as is well-known. A similar remark applies to the $A$-structure of the variation. (It is possible, however, to avoid the use of the local coordinate for the construction of the $A$-structure in the proof of Theorem~3, if one uses the Lefschetz decomposition (4.5.5) together with the fact that $N$ is globally well-defined on $M_{D'}$. Indeed, we have
$$\Gr^W_{w-k}({\rm Ker}\,N)=N^k(P_N\Gr^W_{w+k}M_{D'})\subset\Gr^W_{w-k}M_{D'},$$
and there is a well-defined $A$-structure on the left-hand side given by $\Gr^W_{w-k}(i'{}^*j'_*L)$, where $i:D'\into Y$ and $j':Y\setminus E\into Y$ are natural inclusions. It seems that we have to use the coordinate $y$ for the proof of the compatibility between the dualities of local systems and filtered vector bundles via (4.5.1).)
\ms
(iii) It is possible to prove Theorem~3 without showing that the $A$-structure and the polarization are independent of the choice of a local coordinate if one uses a Zariski-local defining function of the divisor $\rho(D)\subset Y$. Indeed, take a function on a Zariski-open subset $V$ of $Y$ which defines $\rho(D)\cap V$ in $V$. We then get an $A$-structure and a polarization on $\Gr^W_k(M,F)|_{D'\cap V}$ by the same argument as in the proof. We can show that they naturally extend over $D'$ by using [Sch] for the non-degenerate degeneration case (i.e., $N_i=0$ for any $i$). Here the $A$-local system extends since the associated $\C$-local system does by the compatibility with the induced connection on $\Gr^W_kM_{D'}$.
(It is also possible to replace $(Y,E)$ by using Remark~(4.4)(ii).)
\ms
(iv) There is a notion of a mixed Hodge module of geometric origin. This can be obtained by repeating the standard cohomological functors $H^if_*$, $H^if_!$, $H^if^*$, $H^if^!$, $\D$, $\boxtimes$, etc., to $A_{h,pt}$, and taking direct sums and subquotients in the category of mixed Hodge modules.
One can show that each graded piece of a mixed Hodge module of geometric origin on $Y$ is isomorphic to a direct sum of direct factors of certain cohomological direct images of the constant sheaves $A_{h,X_j}$ by projective morphisms $f_j:X_j\to Y$ where the $X_j$ are smooth. (This follows, for instance, from [Sa6], Prop.~7.2 together with the calculation of the vanishing cycles with unipotent monodromy in the normal crossing case [Sa1] which gives the direct factors appearing in the decomposition theorem.) Using this, one can show that the restriction to a curve of a pure Hodge module of geometric origin is also written in a similar way. (Here it seems also possible to proceed by induction on $\dim Y$ assuming that the inverse image of a Zariski-open subset of each irreducible component of $E\subset Y$ is a relatively simple normal crossing divisor over the image.) This may be used to reduce the semi-positivity theorem for admissible variations of mixed Hodge structure of geometric origin to [Ft].
(Note that we can decrease the degree $i$ of the higher direct image $R^if_*\omega_{X/Y}$ as in [Ko2, Cor.~2.24].)
\ms
(v) It is easy to see that Theorem~3 does not hold without the assumption on the unipotence of the local monodromies as is remarked in [FF]. (In fact, let $Y$ be a $\PP^1$-bundle over $\PP^1$ having two disjoint sections with self-intersection numbers $2$ and $-2$. Take the direct image of the constant sheaf by a double covering $f:X\to Y$ which is ramified only along the two sections.)
\ms\nin
{\bf 4.7.~Remark.} The main theorems in this paper are also valid in the analytic case where $f:X\to Y$ is a proper K\"ahler morphism, $Y$ is a complex manifold, and $(X,D)$ is a reduced simple normal crossing analytic pair. The latter means that $(X,D)$ is locally isomorphic to $(X',X'\cap D')$ with $X',D'$ reduced divisors on a complex manifold having no common irreducible components and such that their union is a normal crossing divisor, and moreover the global irreducible components of $X,D$ are smooth. As for the K\"ahler condition, it is enough to assume that each global irreducible component of $X$ is K\"ahler (or more generally, for each global irreducible component of $X$, there is a second cohomology class with $\R$-coefficients which is represented locally over $Y$ by a K\"ahler form).
\sk
Since Theorem~1 concerns only the dualizing sheaf and the first nonzero piece of the Hodge filtration of a mixed Hodge module on $X$, the assertion is essentially local on $X$ by using the canonical isomorphism between the two sheaves defined on the smooth part of $X$. (In fact, this generic isomorphism is uniquely extended to an isomorphism of the whole sheaves by the calculation of extension groups in (1.3).) Then Theorems~2 and 3 and Corollaries~1 and 2 hold where $Y$ is assumed compact in Theorem~3 and Corollary~2. We use [Sa3, Thm.~0.5] for the direct images of the constant sheaves by proper K\"ahler morphisms. (Note that its proof does not use (0.10) in \lc)

\end{document}